\documentclass[11pt]{article}

\usepackage{amscd,amsmath, amssymb, fancyhdr, url,epsfig, color, enumitem, tocloft}
\usepackage{dutchcal}
\usepackage{graphicx}

\usepackage[backref=page]{hyperref}
\renewcommand*{\backrefalt}[4]{%
	\ifcase #1 (Not cited.)%
	\or        (Cited on page~#2.)%
	\else      (Cited on pages~#2.)%
	\fi}

\hypersetup{
	colorlinks   = true,
	citecolor    = magenta,
	linkcolor    =blue,
	urlcolor     =magenta	
}

\newcommand{\version}{version 1.5,\ \ March 4, 2020}
\makeatletter
\def\x@arrow{\DOTSB\Relbar}
\def\xlongrightarrowfill@{\arrowfill@\relbar\relbar\longrightarrow}
\newcommand{\xlongrightarrow}[2][]{%
        \ext@arrow 0099\xlongrightarrowfill@{#1}{#2}}
\makeatother

\newcommand{\la}{\lambda}

\newcommand{\be}{\beta}

\newcommand{\ka}{K\"ahler}
\newcommand{\f}{\varphi}

\newcommand{\Lie}{\operatorname{Lie}}

\newcommand{\vol}{\operatorname{vol}}

\numberwithin{equation}{section}

\def\eqref#1{(\ref{#1})}
\newcommand{\goth}{\mathfrak}
\newcommand{\g}{{\mathfrak g}}

\newcommand{\C}{{\mathbb C}}
\newcommand{\R}{{\mathbb R}}

\newcommand{\6}{\partial}
\def\1{\sqrt{-1}\:}

\newcommand{\cntrct}                % contraction with a vector field
{\hspace{2pt}\raisebox{1pt}{\text{$\lrcorner$}}\hspace{2pt}}
\newcommand{\arrow}{{\:\longrightarrow\:}}

\renewcommand{\bar}{\overline}
\renewcommand{\phi}{\varphi}
\renewcommand{\epsilon}{\varepsilon}
\renewcommand{\geq}{\geqslant}
\renewcommand{\leq}{\leqslant}

\newcommand{\im}{\operatorname{im}}
\newcommand{\End}{\operatorname{End}}

\newcommand{\const}{\operatorname{\text{\sf const}}}

\newcommand{\Iso}{\operatorname{Iso}}
\newcommand{\Aut}{\operatorname{Aut}}

\newcommand{\Tw}{\operatorname{Tw}}

\newcounter{Mycounter}[section]
\newcounter{lemma}[section]
\newcounter{claim}[section]
\newcounter{sublemma}[section]
\newcounter{corollary}[section]
\renewcommand{\thecorollary}{{Corollary \thesection.\arabic{corollary}}}
\newcommand{\corollary}{%
     \setcounter{corollary}{\value{Mycounter}}
     \refstepcounter{corollary}
     \stepcounter{Mycounter}
     {\noindent \bf \thecorollary:\ }}

\newcounter{theorem}[section]
\renewcommand{\thetheorem}{{Theorem \thesection.\arabic{theorem}}}
\newcommand{\theorem}{%
     \setcounter{theorem}{\value{Mycounter}}
     \refstepcounter{theorem}
     \stepcounter{Mycounter}
     {\noindent \bf \thetheorem:\ }}

\newcounter{conjecture}[section]
\renewcommand{\theconjecture}{{Conjecture \thesection.\arabic{conjecture}}}
\newcommand{\conjecture}{%
     \setcounter{conjecture}{\value{Mycounter}}
     \refstepcounter{conjecture}
     \stepcounter{Mycounter}
     {\noindent \bf \theconjecture:\ }}

\newcounter{proposition}[section]
\renewcommand{\theproposition} {{Proposition \thesection.\arabic{proposition}}}
\newcommand{\proposition}{%
     \setcounter{proposition}{\value{Mycounter}}
     \refstepcounter{proposition}
     \stepcounter{Mycounter}
     {\noindent \bf \theproposition:\ }}

\newcounter{definition}[section]
\renewcommand{\thedefinition} {{Definition~\thesection.\arabic{definition}}}
\newcommand{\definition}{%
     \setcounter{definition}{\value{Mycounter}}
     \refstepcounter{definition}
     \stepcounter{Mycounter}
     {\noindent \bf \thedefinition:\ }}

\newcounter{example}[section]
\renewcommand{\theexample}{{Example \thesection.\arabic{example}}}
\newcommand{\example}{%
     \setcounter{example}{\value{Mycounter}}
     \refstepcounter{example}
     \stepcounter{Mycounter}
     {\noindent \bf \theexample:\ }}

\newcounter{remark}[section]
\renewcommand{\theremark}{{Remark \thesection.\arabic{remark}}}
\newcommand{\remark}{%
     \setcounter{remark}{\value{Mycounter}}
     \refstepcounter{remark}
     \stepcounter{Mycounter}
     {\noindent \bf \theremark:\ }}

\newcounter{problem}[section]
\newcounter{question}[section]
\makeatletter

\@addtoreset{equation}{section}
\@addtoreset{footnote}{section}
\makeatother

\def\blacksquare{\hbox{\vrule width 5pt height 5pt depth 0pt}}
\def\endproof{\blacksquare}

\newcommand{\proof}{{\bf Proof: \ }}
\newcommand{\pstep}{{\bf Proof. Step 1: \ }}

\begin{document}

\begin{center}
{\Large\bf  Twisted Dolbeault cohomology of nilpotent Lie algebras}\\[5mm]

%%%%%%%%%%%%%%%%%%%%%%%%%%%%%%%%%%%%%%%%%%%%%%%%%%%%%%%%%%%%%%%%%%%%%%%%%%%%%%%%%%%%
{\large
Liviu Ornea\footnote{Liviu Ornea is  partially supported by a grant of Ministry of Research and Innovation, CNCS - UEFISCDI, project number PN-III-P4-ID-PCE-2016-0065, within PNCDI III.},  
Misha
Verbitsky\footnote{Misha Verbitsky is partially supported by the 
Russian Academic Excellence Project '5-100'', FAPERJ E-26/202.912/2018 
and CNPq - Process 313608/2017-2.\\[1mm]
\noindent{\bf Keywords:} nilpotent Lie algebra, twisted cohomology, locally conformally K\"ahler, Vaisman manifold.

\noindent {\bf 2010 Mathematics Subject Classification:} {53C55, 18H25, 22E25.}
}\\[4mm]

}

\end{center}

{\small
\hspace{0.15\linewidth}
\begin{minipage}[t]{0.72\linewidth}
{\bf Abstract} \\ 
It is well known that the cohomology of any non-trivial 1-dimensional
local system on a nilmanifold vanishes (this result, due to J. Dixmier, was also announced and proved in some particular case by Alaniya). 
A complex nilmanifold is a quotient of a nilpotent Lie group equipped 
with a left-invariant complex structure by an action of a discrete, 
co-compact subgroup. We prove a Dolbeault version 
of Dixmier's and Alaniya's theorem, showing that the Dolbeault cohomology $H^{0,p}(\goth{g},L)$ 
of a nilpotent Lie algebra with coefficients in any non-trivial 1-dimensional
local system vanishes. Note that the Dolbeault cohomology of the
corresponding local system on the manifold is not necessarily zero.
This implies that the twisted version of Console-Fino theorem is false
(Console-Fino proved that the Dolbeault cohomology of a complex nilmanifold
is equal to the Dolbeault cohomology of its Lie algebra, when the complex structure is rational).
As an application, we give a new
proof of a theorem due to H. Sawai, who obtained an explicit
description of LCK nilmanifolds. An LCK structure on a manifold
$M$ is a  K\"ahler structure on its cover $\tilde M$ such that the
deck transform map acts on $\tilde M$ by homotheties.
We show that any complex nilmanifold admitting an LCK structure
is Vaisman, and is obtained as a compact quotient of the product of a Heisenberg 
group and the real line.
\end{minipage}
}
%%%%%%%%%%%%%%%%%%%%%%%%%%%%%%%%%%%%%%%%%%%%%%%%

\tableofcontents

%%%%%%%%%%%%%%%%%%%%%%%%%%%%%%%%%%%%%%%%%%%%%%%%%%%

\section{Introduction}
\label{_Intro_Section_}
%%%%%%%%%%%%%%%%%%%%%%%%%%%%%%%%%%%%%%%%%%%%%%%%%%%

In complex dimension $>2$, an LCK (locally conformally K\"ahler) 
manifold is a complex manifold equipped with a Hermitian metric 
(LCK metric) which is locally conformally equivalent to a K\"ahler manifold.
Then the Hermitian form satisfies $d(\omega)=\theta \wedge \omega$, where
$\theta$ is a 1-form. Clearly,  $d^2(\omega)=0$ gives $d\theta \wedge \omega=0$,
which is equivalent to $d\theta=0$ in complex dimension $>2$.
In complex dimension 2, this condition has to be added artificially.

Summing it up, {\bf an LCK metric} (or {\bf an LCK structure}) 
on a complex manifold of complex
dimension $\geq 2$ is a Hermitian metric with Hermitian form
$\omega$ which satisfies $d(\omega)=\theta \wedge \omega$, where
$\theta$ is a closed 1-form. 

Consider the differential $d_\theta(\alpha):= d(\alpha) - \theta \wedge \alpha$.
Then the equation $d(\omega)=\theta \wedge \omega$ can be written as $d_\theta(\omega)=0$.
The differential $d_\theta$ is called {\bf the Morse-Novikov differential}.
It can be interpreted as de Rham differential for forms with coefficients
in a local system, which is done as follows.

Let $L$ be a trivial real bundle of rank 1, and $\nabla_0$
the standard connection. Demote by $\nabla:=\nabla_0-\theta$ the connection
$\nabla(b)=\nabla_0(b)- b\otimes \theta$. Since $\nabla^2=0$,
this connection is flat, and $(L, \nabla)$ defines a local system,
which we denote by the same letter. Then $d_\theta$ is the de Rham differential
on differential forms with values in $L$, and the cohomology of $d_\theta$
is identified with the cohomology of this local system. 

In this paper we study  the locally conformally K\"ahler structures on nilmanifolds.

\hfill

A nilmanifold is a quotient of a connected, simply connected 
nilpotent Lie group by a discrete co-compact lattice. It is not
hard to see that any nilmanifold is the total space of an 
iterated family of circle bundles. Its cohomology
is expressed through the cohomology of the Chevalley-Eilenberg differential
on the corresponding Lie algebra $\g$, and the Morse-Novikov cohomology is  obtained
as cohomology of the twisted version of the Chevalley-Eilenberg differential
(Section \ref{_Twisted_cohomologies_}).

In the literature, there are several different (and incompatible)
ways to define a complex nilmanifold. Before 1980-ies, most
people defined complex nilmanifold as a  quotient of a complex nilpotent Lie group
$G_\C$ by a co-compact lattice. This quotient is a parallelizable
complex manifold (that is, its tangent bundle is holomorphically
trivial). Indeed, the right action of $G_\C$  on itself
commutes with the left action, and both are holomorphic.
If we take the quotient $M:=G_\C/\Gamma$ by (say) left action
of $\Gamma\subset G_\C$, the manifold $M$ remains homogeneous
with respect to the right action of $G_\C$, and this action
trivializes $TM$, making $M$ parallelizable.

Now such nilmanifolds are called {\bf parallelizable complex nilmanifolds}
or {\bf Iwasawa type nilmanifolds}, after the Iwasawa group (the non-commutative 
3-dimensional complex nilpotent Lie group).

Before the paper \cite{_Cordeiro_Fernandez_Deleon_}, in the published literature
``complex nilmanifold'' meant the quotient of a complex Lie group;
see for example \cite{_Fisher_} or \cite{_Oeljeklaus:hyperflachen_}.
In \cite{_Cordeiro_Fernandez_Deleon_}, Cordero, Fern\'andez
and de Le\'on introduced the modern notion of a complex
nilmanifold. Since then, a complex nilmanifold is defined as a compact quotient of a nilpotent
Lie group equipped with a left-invariant complex 
structure (Subsection \ref{_LCK_nilma_Subsection_}). The main
advantage of this notion is that it can be rephrased
in terms of a  complex structure operator on the corresponding
Lie algebra.
However, such a quotient is not homogeneous, unless
$G$ is a complex Lie group.

The same approach can be used to describe other 
important geometric structures on nilmanifolds in terms
of linear-algebraic structures on their Lie algebras 
(Subsection \ref{_inva_Lie_Subsection_}).
In the same paper \cite{_Cordeiro_Fernandez_Deleon_},
this approach was used to construct an
LCK structure on a nilmanifold obtained
from the Heisenberg group. Later, H. Sawai
proved in \cite{saw1} that any LCK structure on a complex
nilmanifold is obtained this way. In the present
paper we give a proof of Sawai's theorem
based on cohomology vanishing.

Consider an oriented rank 1 local system 
$(L, \nabla_0+\theta)$ on a nilmanifold $M=G/\Gamma$, where 
$(L, \nabla_0)$ is a trivial bundle with connection, and
$\theta$ is
a closed, non-exact 1-form. In \cite{_Alaniya_} (see also
\cite{mili}), it was shown that $H^*(M, L)=0$.
The proof goes as follows: first, one identifies 
$H^*(M, L)$ with the cohomology of the Lie algebra
$\g=\Lie(G)$ with
coefficients in a non-trivial rank one representation.
This is done by using the homogeneous $G$-action.
Then one takes the Chevalley-Eilenberg complex
which computes the Lie algebra cohomology, 
filters it by the central series, and 
computes the first page of the corresponding
spectral sequence. This first page
is identified with the twisted cohomology of a commutative
Lie algebra, which always vanishes. 

However, this approach will not work for Dolbeault
cohomology.  Console and Fino (\cite{_Console_Fino_})
proved that the Dolbeault cohomology can, indeed, be
computed using the Hodge decomposition on the
Chevalley-Eilenberg complex of the
corresponding Lie algebra, when the complex structure is rational (in dimension up to 6 the rationality of the complex structure is not needed, see \cite{_Fino_Rollenske_Ruppenthal_}, but for higher dimensions the question is still open). However, this result is
highly non-trivial, because complex nilmanifolds are
not homogeneous. Moreover, the extension of Console-Fino theorem
to the Dolbeault cohomology with coefficients 
in a local system (``twisted Dolbeault cohomology'') is false
(Subsection \ref{_Console-Fino_coeff_Subsection_}).

The main result of the present paper is the following theorem,
which computes the Lie algebra version of the Dolbeault cohomology
with coefficients in a local system.

\hfill

%%%%%%%%%%%%%%%%%%%%%%%%%%%%%%%%%%%%%%%%%%%%%%%%%%%%%%%%%%%%%%%%%%%%%%%%
\theorem
Let $\g$ be a nilpotent real Lie algebra, $I:\; \g \arrow \g$
a complex structure operator, $\g_\C:=\g\otimes_\R \C$, and $\g_\C= \g^{1,0}\oplus \g^{0,1}$
the corresponding eigenvalue decomposition, called the
Hodge decomposition in the sequel. Assume that
$\g^{1,0}\subset \g_\C$ is a Lie subalgebra (this is
equivalent to the left-invariant complex structure
on on the Lie group $G=\Lie(\g)$ associated with $I$ 
being integrable, see \ref{_complex_on_nil_alg_Definition_}).
Consider the Hodge decomposition $\Lambda^*(\g_\C^*)= \bigoplus\Lambda^{p,q}(\g^*)$ on the
Grassman algebra of $\g_\C$, with  
\[
\Lambda^{p,q}(\g^*)=\Lambda^p((\g^{1,0})^*)\otimes_\C \Lambda^q((\g^{0,1})^*)
\]
and define the twisted Dolbeault differentials $\6_\theta$, $\bar\6_\theta$ as Hodge components of the
twisted Chevalley-Eilenberg differential $d_\theta(x):= d(x)- \theta \wedge x$, i.e. $\6_\theta=\6-\theta^{1,0}$ and $\bar\6_\theta=\bar\6-\theta^{0,1}$ 
(Subsection \ref{_inva_Lie_Subsection_}).
Then the cohomology of the complex $(\Lambda^{0,*}(\g_\C^*), \bar\6_\theta)$ vanishes.

\proof See \ref{_twi_Dolbeault_Vanish_}. \endproof

\hfill

We apply this result to get a new proof of the classification
of LCK structures on complex nilmanifolds.

%%%%%%%%%%%%%%%%%%%%%%%%%%%%%%%%%%%%%%%%%%%%%%%%%%%%%%%%%%%%%%%%%%%%%%%%%%%%%%%%%%%%

\section{Preliminaries}

%%%%%%%%%%%%%%%%%%%%%%%%%%%%%%%%%%%%%%%%%%%%%%%%%%%%%%%%%%%%%%%%%%%%%%%%%%%%%%%%%%%%

%%%%%%%%%%%%%%%%%%%%%%%%%%%%%%%%%%%%%%%%%%%%%%%%%%%%%%%%%%%%%%%%%%%%%%%%
\subsection{Locally conformally K\"ahler manifolds}
%%%%%%%%%%%%%%%%%%%%%%%%%%%%%%%%%%%%%%%%%%%%%%%%%%%%%%%%%%%%%%%%%%%%%%%%

\definition Let $(M,I)$ be a complex manifold, $\dim_\C M\geq 2$. It
is called {\bf locally conformally K\"ahler} (LCK) if it
admits a  Hermitian metric $g$ whose 
fundamental 2-form $\omega(\cdot,\cdot):=g(\cdot, I\cdot)$
satisfies
\begin{equation}\label{deflck}
d\omega=\theta\wedge\omega,\quad d\theta=0,
\end{equation}
for a certain closed 1-form $\theta$ called {\bf the Lee
	form}. 

\hfill

\remark 
Definition \eqref{deflck} is equivalent to the existence
of a  covering $\tilde M$ endowed with a K\"ahler metric $\Omega$ which is
acted on by the deck group $\Aut_M(\tilde M)$ by homotheties. Let 
\begin{equation}\label{chi}
\chi:\Aut_M(\tilde M)\longrightarrow \R^{>0},\quad \chi(\tau)=\frac{\tau^*\Omega}{\Omega},
\end{equation}
be the character which associates to a homothety its scale factor. On this cover, the pull-back of the Lee form is exact.

\hfill

\remark\label{_minimal_cover_}  For an LCK manifold, coverings with the above property are not unique. The covering for which the character $\chi$ is injective is called the {\bf minimal cover}.

\hfill

\remark The operator $d_\theta:=d-\theta\wedge$ obviously satisfies $d_\theta^2=0$ and hence $(\Lambda^*M,d_\theta)$ produces a cohomology called {\bf Morse-Novikov} or {\bf twisted}. It can be interpreted as the cohomology of the local system $L$ associated to the  line bundle endowed with (flat) connection form $\theta$. See Section \ref{_Twisted_cohomologies_} for details.

\hfill

\theorem (\cite{va_tr})\label{vailcknotk}
Let $(M,\omega, \theta)$ be a compact LCK manifold, not globally conformally K\"ahler (i.e. with non-exact Lee form). 
Then $M$ does not admit a K\"ahler metric.

%%%%%%%%%%%%%%%%%%%%%%%%%%%%%%%%%%%%%%%%%%%%%%%%%%%%%%%%%%%%%%%%%%%%%%%%
\subsection{Vaisman manifolds} 
%%%%%%%%%%%%%%%%%%%%%%%%%%%%%%%%%%%%%%%%%%%%%%%%%%%%%%%%%%%%%%%%%%%%%%%%

\definition An LCK manifold $(M,\omega, \theta)$ is called {\bf
	Vaisman} if $\nabla\theta=0$, where $\nabla$ is the
Levi-Civita connection of $g$. 

\hfill

The following characterization, very much used in applications, is available:

\hfill

\theorem (\cite{kor}) Let $(M,\omega, \theta)$ be an LCK manifold equipped with a 
holomorphic and conformal $\C$-action  without fixed points,
which lifts to non-isometric homotheties on 
the K\"ahler covering $\widetilde M$. { Then $(M,\omega, \theta)$
	is conformally equivalent to a Vaisman manifold.}

\hfill

The main example of Vaisman manifold is the
diagonal Hopf manifold (\cite{_OV:Shells_}). The Vaisman compact complex
surfaces are classified in \cite{_Belgun_}, see also \cite{_ovv_surfaces_}. 

\hfill

\remark There exist compact LCK manifolds which do not admit Vaisman metrics. Such are the LCK Inoue surfaces, \cite{_Belgun_}, the Oeljeklaus-Toma manifolds, \cite{_Oeljeklaus_Toma_}, \cite{_Otiman_}, and the non-diagonal Hopf manifolds, \cite{_OV:Shells_}, \cite{_ovv_surfaces_}.

\hfill

\remark\label{_canon_foli_totally_geodesic_Remark_}
On a Vaisman manifold, the Lee field $\theta^\sharp$ and the anti-Lee field $I\theta^\sharp$ are real holomorphic ($\Lie_{\theta^\sharp}I=\Lie_{I\theta^\sharp}I=0$) and Killing ($\Lie_{\theta^\sharp}g=\Lie_{I\theta^\sharp}g=0$), see \cite{_Dragomir_Ornea_}. 

Recall that a Killing vector field $X$  satisfies the equation: 
\[ g(\nabla_AX, B)= - g(A, \nabla_BX).\]
Setting 
$A=X$, we get $g(\nabla_XX, B)= - g(X, \nabla_BX)$.
The last term is equal to $-1/2 \Lie_B(g(X,X))=0$,
hence $g(\nabla_XX, B)= 0$ for all $B$.
Then $\nabla_XX= 0$, and hence the 
trajectories of a Killing field of constant length are
geodesics. Therefore, the canonical foliation in
a Vaisman manifold is totally geodesic.

\hfill

\remark {\bf (i)} Note that while the LCK condition is conformally invariant (changing the metric $g\mapsto e^fg$ changes the Lee form into $\theta\mapsto \theta+df$), the Vaisman condition is not conformally invariant. Indeed, on a Vaisman manifold, the Lee form is coclosed (being parallel) and hence  a Vaisman metric is a Gauduchon metric; but one knows that on a compact complex manifold, each conformal class contains a unique Gauduchon metric (up to constant multipliers), \cite{_Gauduchon_1984_}.

{\bf (ii)} Since $\theta$ is parallel, it has constant norm and thus we can always scale the LCK metric such that $|\theta|=1.$ In this assumption, the following formula holds, \cite{_Vaisman:gen_hopf_}, \cite{_Dragomir_Ornea_}:
\begin{equation}\label{dctheta}
d\theta^c=\theta\wedge\theta^c-\omega, \quad \text{where}\quad \theta^c(X)=-\theta(IX).
\end{equation}
Moreover, one can see, \cite{_Verbitsky:LCHK_}, that the
(1,1)-form $\omega_0:=-d^c\theta$ is semi-positive
definite, having all eigenvalues\footnote{The eigenvalues
	of a Hermitian form $\eta$  are the eigenvalues of the
	symmetric operator $L_\eta$ defined by the equation
	$\eta(x,Iy)=g(L_\eta x,y)$.} 
positive, except one which is 0.

\subsection{LCK manifolds with potential}

We now introduce a class of LCK manifolds strictly containing  the Vaisman manifolds.

\hfill

\definition We say that an LCK manifold has  {\bf LCK potential} if it
admits a K\"ahler covering on which the K\"ahler metric
has a global and positive  potential function which
is acted on by holomorphic homotheties by the deck group. In this case, $M$ is called {\bf LCK manifold with potential}.

\hfill

\remark Note that in several previous papers of ours we asked the potential function to be proper. It is not the case for the above definition. See \cite{ov_jgp_16} for a comprehensive discussion about proper and improper LCK potentials.

\hfill

\remark\label{dc_on_pot} One can prove that $(M,I,g,\theta)$ is LCK with potential if and only if equation \eqref{dctheta} is satisfied.  

\hfill

\definition A function $\phi\in C^\infty (M)$ is called {\bf $d_\theta d^c_\theta$-plurisubharmonic} if $\omega= d_\theta d^c_\theta(\phi)$.

\hfill

\remark Note that $d_\theta d^c_\theta$-plurisubharmonic are not plurisubharmonic (they exist on compact manifolds). Moreover, if $\phi$ is $d_\theta d^c_\theta$-plurisubharmonic, then $\phi+\const$ is not necessarily $d_\theta d^c_\theta$-plurisubharmonic.

\hfill

Equation \eqref{dctheta} (and hence the definition of LCK manifolds with potential) can be translated on the LCK  manifold itself:

\hfill

\theorem (\cite[Claim 2.8]{{ov_jgp_16}})\label{_theta_pluri_iff_pluri_} $(M,I,\theta,\omega)$ is LCK with potential if and only if $\omega= d_\theta d^c_\theta(\phi)$ for a strictly positive $d_\theta d^c_\theta$-plurisubharmonic function $\phi$ on $M$. 

\hfill

\theorem (\cite{ov_poz_pot}) \label{_strictly_negative_pots_Theorem_}
Let $(M, \theta, \omega)$ be an LCK manifold
which is not K\"ahler,
and suppose $\omega= d_\theta d^c_\theta(\phi_0)$ for some
smooth   function
$\phi_0\in C^\infty (M)$.
Then $\phi_0>0$ at some point of $M$.

\hfill

\proof (Courtesy of Matei Toma.) By absurd,  $\f_0\leq 0$
everywhere on $M$. Let $\tilde M$ be the minimal 
K\"ahler cover of $M$ (see \ref{_minimal_cover_}), $\Gamma$ its deck group, and $\rho$ a function
on $\tilde M$ such that $d\rho=\theta$.
Then the $d_\theta d^c_\theta$-plurisubharmonicity of
$\phi_0$ is equivalent to the plurisubharmonicity of
$\phi:=e^{-\rho}\phi_0$, see \ref{_theta_pluri_iff_pluri_}.

Since the strict $d_\theta d^c_\theta$-plurisubharmonicity is stable under $C^2$-small deformations
of $\phi$, the function $\phi-\epsilon$ is also strictly
$d_\theta d^c_\theta$-plurisubharmonic. Therefore, we may 
assume that $\phi <0$ everywhere.

Define
$$\psi:=-\log(-\f).$$ 
Since $x\arrow -\log(-x)$ is strictly
monotonous and convex, the function
$\psi$ is strictly $d_\theta d^c_\theta$-plurisubharmonic.  Moreover,
for every element $\gamma\in\Gamma$, we have
\[
\gamma^*\psi=-\log(-(\f\circ\gamma))=-
\log(\chi(\gamma))-\log(-\f)=\const+\psi.
\]
Therefore, the K\"ahler form
$dd^c\psi$ is $\Gamma$-invariant and descends to $M$, contradiction with \ref{vailcknotk}.
\endproof

\hfill

\remark All Vaisman manifolds are LCK manifolds with potential. Among the non-Vaisman examples, we mention the non-diagonal Hopf manifolds, \cite{ov_jgp_16}. On the other hand, the Inoue surfaces  and their higher dimensional analogues, the Oeljeklaus-Toma manifolds, are   compact LCK manifolds which are not LCK manifolds with potential  (\cite{_Otiman_}).

\hfill

We can characterize the Vaisman metrics among the LCK metrics with potential:

\hfill

\proposition\label{pot_gau}(\cite[Proposition 2.3 \& Corollary 2.4]{ov_hopf_surf})\\ 
Let $(M,\omega,\theta)$ be a compact LCK manifold with
potential. Then the LCK metric is Gauduchon if and only if
$\omega_0=-d^c\theta$ is semi-positive definite, and hence it is Vaisman. 
Equivalently, a compact LCK manifold with potential and with
constant norm of $\theta$ is Vaisman.

%%%%%%%%%%%%%%%%%%%%%%%%%%%%%%%%%%%%%%%%%%%%%%%%%%%%%%%%%%%%
\section{LCK structures on nilmanifolds}

\subsection{Invariant geometric structures on Lie groups}
\label{_inva_Lie_Subsection_}
%%%%%%%%%%%%%%%%%%%%%%%%%%%%%%%%%%%%%%%%%%%%%%%%%%%%%%%%%%%%%%%%%%%%%%%%

Let $G$ be a Lie group, $\Lambda$ a discrete subgroup, and
$I$ a left-invariant complex structure on $G$. Consider
the quotient space $G/\Lambda$, where $\Lambda$ acts by
left translations. Since $I$ is left-invariant,
the manifold $G/\Lambda$ is equipped with a natural
complex structure. This construction is usually 
applied to solvmanifolds or nilmanifolds, but
it makes sense with any Lie group. When we need to refer to this
particular kind of complex structures on $G/\Lambda$, we call them 
{\bf locally $G$-invariant}.

\hfill

A left-invariant almost complex structure $I$ on $G$ is
determined by its restriction to the Lie algebra
$\g=T_e G$, giving the Hodge decomposition
$\g\otimes_\R \C= \g^{1,0}\oplus \g^{0,1}$. 
This structure is integrable
if and only if the commutator of $(1,0)$-vector fields
is again a $(1,0)$-vector field, which is equivalent
to $[\g^{1,0},\g^{1,0}]\subset \g^{1,0}$.
This allows one to define a complex structure
on a Lie algebra.

\hfill

\definition\label{_complex_on_nil_alg_Definition_}
A complex structure on a Lie algebra $\g$ is a subalgebra
$\g^{1,0}\subset \g\otimes_\R \C$
which satisfies $\g^{1,0}\oplus
\overline{\g^{1,0}}=\g\otimes_\R \C$.

\hfill

Indeed, an almost complex structure operator
$I$ can be reconstructed from the decomposition
$\g^{1,0}\oplus\overline{\g^{1,0}}=\g\otimes_\R \C$ by making it act
as $\1$ on $\g^{1,0}$ and $-\1$ on $\overline{\g^{1,0}}$.

\hfill

In a similar way one could define 
symplectic structures or LCK structures on a Lie algebra $\g$.
Recall that the Grassmann algebra\index{Grassmann algebra} $\Lambda^*(\g^*)$ is
equipped with a natural differential, called 
{\bf Chevalley differential}, which is equal 
to the de Rham differential on left-invariant
differential forms on the Lie group if we identify
the space of such forms with $\Lambda^*(\g^*)$. 
In dimensions 1 and 2 it can be written explicitly as follows:
if $\la\in\goth{g^*}$, then $d\la(x,y)=-\la [x,y]$; if
$\be\in\Lambda^2\goth{g}^*$, then $d\be(x,y,z)=-\be([x,y
], z) - \be([y,z], x) - \be([z,x], y )$.
The corresponding complex is called 
{\bf the Chevalley-Eilenberg complex}.

Further on, we shall always interpret the elements
of $\Lambda^p(\g^*)$ as left-invariant differential
forms on the corresponding Lie group, and refer
to them as to ``$p$-forms'', with all the usual
terminology (``closed forms'', ``exact forms'')
as used for the elements of de Rham algebra.

\hfill

\definition
{\bf A symplectic structure} on a Lie algebra $\g$ is a
non-degenerate, closed 2-form $\omega\in
\Lambda^2(\g^*)$. A {\bf K\"ahler structure}
on a Lie algebra $\g$ is a complex structure $I$
and a Hermitian form $h$ on $\g$ such that 
the fundamental 2-form 
$\omega(\cdot, \cdot):= h(\cdot, I(\cdot))$
is closed. 

\hfill

%%%%%%%%%%%%%%%%%%%%%%%%%%%%%%%%%%%%%%%%%%%%%%%%%%%%%%%%%%%%
\remark
In \cite{bens_gor} 
(see also \cite{hase}) it was shown that any
nilpotent Lie algebra admitting a K\"ahler structure is
actually abelian.\index{Lie algebra!nilpotent!}

\subsection{LCK nilmanifolds}\label{_LCK_nilma_Subsection_}

\definition\label{_LCK_Lie_algebra_Definition_}
{\bf An LCK structure} on a Lie algebra is a
complex structure $I$
and a Hermitian form $h$ on $\g$ such that 
the fundamental 2-form 
$\omega(\cdot, \cdot):= h(\cdot, I(\cdot))$
satisfies $d(\omega) = \theta \wedge \omega$, where
$\theta\in \Lambda^1(\g^*)$ is a closed 1-form.

\hfill

\definition
A {\bf nilmanifold}\index{nilmanifold} is a 
quotient of a simply connected nilpotent\index{Lie group!nilpotent}
Lie group by a co-compact discrete
subgroup. Alternatively (\cite{mal}), one can define
nilmanifolds as  manifolds which
admit a homogeneous action\index{action!homogeneous} by a nilpotent
Lie group.

\hfill

%%%%%%%%%%%%%%%%%%%%%%%%%%%%%%%%%%%%%%%%%%%%%%%%
\definition\label{_LCK_nilma_Definition_}
Let $(G, I)$ be a nilpotent %or solvable 
Lie group\index{Lie group!solvable}\index{Lie group!nilpotent}
equipped with a left-invariant complex structure.
For any co-compact discrete
subgroup $\Gamma\subset G$, the (left) quotient
$G/\Gamma$ is equipped with a natural complex structure.
This quotient is called  {\bf a complex nilmanifold}.
	%(solvmanifold)}. 
Similarly, if $G$ is a Lie group
with left-invariant LCK structure, the quotient
$G/\Gamma$ is called {\bf an LCK nilmanifold}. 
	%(solvmanifold)}. 
	This structure is clearly
locally homogeneous. \index{manifold!locally homogeneous}When we need to refer to this
particular kind of structures, we call them 
{\bf locally $G$-invariant}.

\hfill

\remark The existence of a co-compact lattice already
imposes strong restrictions: the group should be
unimodular\index{Lie group!unimodular},
\cite[\S 1.1, Exercise 14b]{_Morris:Ratner_}. All nilpotent Lie groups are unimodular, but this is not sufficient for the existence of a cocompact lattice. According to \cite{mal}, a nilpotent Lie group admits a cocompact lattice if and only if its Lie algebra admits a basis in which the structural constants are all rational.
%, but not all solvable ones are. 

\hfill

In this paper 
we present results concerning LCK structures on
nilmanifolds. A classification of nilpotent Lie algebras which admit LCK structures is given.

\hfill 

\remark  \label{_Malcev_comple_Remark_}
As shown by  Mal\v cev, \cite{mal},  
any nilmanifold is uniquely determined by its
fundamental group, which is a discrete nilpotent torsion-free group, and any such
group uniquely determines a nilmanifold.

\hfill

\remark Note that right translations are not holomorphic
with respect to left invariant complex structures. This
means that  for a lattice $\Lambda$, the complex structure
induced on the manifold $G/\Lambda$ is not necessarily
invariant, and hence $G/\Lambda$ is not a homogeneous LCK
manifold.\index{manifold!LCK!homogeneous}

\hfill

\remark Note that it may happen that a K\"ahler manifold  be diffeomorphic (as real manifold) with a nilmanifold but the K\"ahler structure need not be invariant. For example, let $E$ be  an elliptic curve which is  a 2-fold ramified covering of  $\C P^1$, and take the fibered product  $E\times_{\C P^1} \Tw(T^4)$, 
where $\Tw(T^4)$ is the twistor space of a 4-torus. This fibered product is
actually a torus, but with an inhomogeneous complex structure
(since it has non-trivial canonical bundle, in fact its anticanonical bundle
is semi-positive, with many non-trivial sections); compare with \cite[\S 5]{_Calabi_}. This example suggests that the same phenomenon could appear for LCK manifolds (to be diffeomorphic with a nilmanifold without having an invariant LCK structure), but for the moment we don't have a concrete example.

\hfill

\example\label{vai_hei}\label{kod_nil}  
Let $H_{2n-1}$ be the Heisenberg group of matrices of the form
$$\begin{pmatrix}
1&A&c\\
0&I_{2n-1}&B^t\\
0&0&1
\end{pmatrix}, \quad c\in\R,\quad A,B\in\R^{n-1}.$$
Its Lie algebra $\goth{h}_{2n-1}$
can be defined in terms of generators and relations as follows:
it has a basis $\{X_i,Y_i,Z\}$, $i=1,\ldots,n-1$.
with $[X_i,Y_i]= Z$, and the rest of Lie brackets trivial.
Clearly, it is a nilpotent Lie algebra. After taking the
quotient of
the corresponding Lie group
by a co-compact discrete subgroup $\Lambda$, one obtains 
a nilmanifold called {\bf the Heisenberg nilmanifold}. It
has an invariant Sasakian structure given by requiring the
above basis to be orthonormal, defining the contact
form as $\sum y_idx_i+dz$ (where $dx_i,dy_i, dz$ are the dual
invariant one-forms). Then $Z$ is its Reeb field.

Now, the product $H_{2n-1}/\Lambda\times S^1$ is a Vaisman
nilmanifold whose universal cover is the product $
H_{2n-1}\times\R$, a Lie group with Lie algebra
$\goth{g}:= \goth{h}_{2n-1}\times\R$, $\R=\langle
T\rangle$,  in which the brackets of the basis
$\{X_i,Y_i,Z, T\}$ are the above, and $T$ is central. 
The linear LCK structure on $\goth{g}$ is given by
asking the basis to be orthonormal, and by defining  the
complex structure $IX_i=Y_i$, $IZ=-T$. The fundamental
form is $\omega=\sum (X_i^*\wedge Y_i^*)-Z^*\wedge T^*.$ A direct
computation shows that:
$d\omega=T^*\wedge \omega,$ and hence the  Lee form is $T^*$
which can be seen directly to be $g$-parallel.

%%%%%%%%%%%%%%%%%%%%%%%%%%%%%%%%%%%%%%%%%%%%%%%%%%%%%%%%%%%%%%%%%%%%%%%%%%%%%%%%%%%%

\section{Twisted Dolbeault cohomology on nilpotent Lie algebras}\label{_Twisted_cohomologies_}

%%%%%%%%%%%%%%%%%%%%%%%%%%%%%%%%%%%%%%%%%%%%%%%%%%%%%%%%%%%%%%%%%%%%%%%%%%%%%%%%%%%

%%%%%%%%%%%%%%%%%%%%%%%%%%%%%%%%%%%%%%%%%%%%%%%%%%%%%%%%%%%%
\subsection{Twisted Dolbeault cohomology}
%%%%%%%%%%%%%%%%%%%%%%%%%%%%%%%%%%%%%%%%%%%%%%%%%%%%%%%%%%%%

Let $\theta$ be a closed 1-form on a complex manifold.
Then $d -\theta$ defines a connection on the trivial bundle.
We denote the corresponding local system by $L$.
The {\bf twisted Dolbeault differentials} $\6_\theta:= (d_\theta)^{1,0}=\6-\theta^{1,0}$ and
$\bar\6_\theta:= (d_\theta)^{0,1}=\bar\6-\theta^{0,1}$
are Morse-Novikov counterparts to the usual Dolbeault differentials.
The cohomology of these differentials corresponds to the 
Dolbeault cohomology with coefficients in the 
holomorphic line bundle $L_\C$ obtained as complexification 
of the local system $L$.

Working with LCK nilmanifolds and nilpotent Lie algebras, 
it is natural to consider their twisted Dolbeault cohomology. 
As shown by Console and Fino (\cite{_Console_Fino_}), the usual 
Dolbeault cohomology of complex nilmanifolds is equal to the 
cohomology of the Dolbeault version of the corresponding 
Chevalley-Eilenberg complex on the Lie algebra. This result
does not hold in the twisted Dolbeault cohomology, as we shall see in \ref{_console_fino_ex_}.

However, as we show in \ref{_twi_Dolbeault_Vanish_} below, the corresponding
Lie algebra cohomology vanishes.

\hfill

\definition
Let $\g$ be a Lie algebra and 
$\theta\in \g^*$ a closed 1-form and \[d_\theta:=d-\theta\in \End(\Lambda^*(\g^*)),\]
where $d$ is the Chevalley differential and $\theta$
denotes the operation of multiplication by $\theta$.
The twisted Dolbeault differentials on
$\Lambda^*(\g^*)$ are $\6_\theta:= (d_\theta)^{1,0}$ and
$\bar\6_\theta:= (d_\theta)^{0,1}$.

\hfill

The cohomology of $d_\theta$ is always zero 
for non-zero $\theta$ (\cite{_Alaniya_,mili}). This theorem
is due to J. Dixmier, \cite{_Dixmier_}. 
The main result of the present section is the following theorem,
which is the Dolbeault version of Dixmier's and Alaniya's theorem.

\hfill

%%%%%%%%%%%%%%%%%%%%%%%%%%%%%%%%%%%%%%%%%%%%%%%%%%%%%%%%%%%%
\theorem\label{_twi_Dolbeault_Vanish_}
Let $\g$ be a nilpotent Lie algebra with a complex structure 
(\ref{_complex_on_nil_alg_Definition_}), and 
$\theta\in \g^*$ a non-zero, closed real 1-form.
Then the cohomology of $(\Lambda^{0,*}(\g^*),\bar\6_\theta)$ vanishes.

\hfill

\proof
Consider the central series of the Lie algebra
$\g^{0,1}$, with $W_0=  \g^{0,1}, W_1=  [W_0, W_0], ..., 
W_k=[W_0, W_{k-1}]$. Let $A_k\subset (\g^{0,1})^*$
be the annihilator of $W_k$. For any 1-form $\lambda\in(\g^{0,1})^*$, one has 
$\bar\6(\lambda)(x,y) = \lambda([x,y])$, hence 
$\bar\6(A_k)\subset \Lambda^2(A_{k-1})$. Consider the filtration
$\Lambda^*(A_1)\subset \Lambda^*(A_2)\subset ...$ on $\Lambda\in(\g^{0,1})^*$.
Since $\bar\6(A_k)\subset \Lambda^2(A_{k-1})$, the operator $\bar\6$
shifts the filtration by 1: $\bar\6(\Lambda^*(A_k))\subset \Lambda^*(A_{k-1})$.

Consider the spectral sequence of the complex $(\Lambda^*(\g^*),\bar\6_\theta)$
filtered by $V_0\subset V_1\subset V_2 \subset ...$, where
$V_k := \Lambda^*(A_k)$.
Since  $\bar\6(V_k) \subset V_{k-1}$, the operator $\bar\6_\theta$
acts on $\bigoplus_k V_k/V_{k-1}$ as multiplication by $\theta^{0,1}$. 
The corresponding associated graded complex, which is $E_0^{*,*}$ of this
spectral sequence, is identified with
$\left(\bigoplus_k V_k/V_{k-1}, \theta^{0,1}\right)$. 
We identify $\bigoplus_k V_k/V_{k-1}$ with the
Grassmann algebra $\bigoplus_k V_k/V_{k-1}=\Lambda^*(\g^*)$.
After this identification, the multiplication
\[ \theta^{0,1}:\;\bigoplus_k V_k/V_{k-1}
\arrow \bigoplus_k V_k/V_{k-1}
\]
becomes multiplication by a 1-form, obtained from $\theta^{0,1}$.
The cohomology of multiplication by a 1-form always vanishes. Then the
$E_1^{*,*}$-page of the spectral sequence vanishes, which
implies vanishing of $H^*(\Lambda^*(\g^{0,*}),\bar\6_\theta)$.
\endproof

\hfill

The following corollary will be used in the classification of LCK
structures on nilpotent Lie algebras. 

\hfill

%%%%%%%%%%%%%%%%%%%%%%%%%%%%%%%%%%%%%%%%%%%%%%%%%%%%%%%%%%%%
\corollary\label{_Hodge_chasing_Corollary_}
Let $(\goth{g},  I)$ be a $2n$-dimensional
nilpotent Lie algebra with complex structure, $\theta\in \Lambda^1(\g^*)$
a non-zero, closed real 1-form, and  $\omega\in \Lambda^{1,1}(\g^*)$
a $d_\theta$-closed (1,1)-form.
Then there exists a 1-form $\tau\in \Lambda^1(\g^*)$
such that  $\omega=d_\theta(\tau)$,
where $d_\theta (I\tau)=0$.

\hfill

\proof
Denote by $H^{1,1}_{d_\theta d^c_\theta}(\g^*)$ the twisted
Bott-Chern cohomology, that is, all $d_\theta$-closed (1,1)-forms
up to the image of $d_\theta d^c_\theta$.
The standard exact sequence
\[ H^{0,1}_{\bar\6_\theta}(\g^*) \oplus H^{1,0}_{\6_\theta}(\g^*)
\stackrel d\arrow H^{1,1}_{d_\theta d^c_\theta}(\g^*) \stackrel \mu \arrow H^2(M)
\]
(\cite{_ov:MN_}, equation (4.7)) implies that
the kernel of $\mu$ vanishes when  
$H^{0,1}_{\bar\6_\theta}(\g^*)=H^{1,0}_{\6_\theta}(\g^*)=0$.
The last equation follows from  \ref{_twi_Dolbeault_Vanish_},
hence $\mu$ is injective. The LCK form $\omega$ belongs to the
kernel of $\mu$, because the $d_\theta$-cohomology of
$\Lambda^*(\g^*)$ vanishes by Diximier and Alanya theorem (\cite{_Alaniya_,mili}).
Therefore, $\omega$ is twisted Bott-Chern exact, i.e.  
$\omega= d_\theta d_\theta^c f$ for some
$f\in \Lambda^0(\g^*)=\R$. Set $\tau:=d^c_\theta f$. Then $I\tau=d_\theta f$ and $d_\theta(I\tau)=0$ as stated. \endproof

%%%%%%%%%%%%%%%%%%%%%%%%%%%%%%%%%%%%%%%%%%%%%%%%%%%%%%%%%%%%
\subsection{Console-Fino theorem with coefficients in a local system fails}
\label{_Console-Fino_coeff_Subsection_}
%%%%%%%%%%%%%%%%%%%%%%%%%%%%%%%%%%%%%%%%%%%%%%%%%%%%%%%%%%%%

%%%%%%%%%%%%%%%%%%%%%%%%%%%%%%%%%%%%%%%%%%%%%%%%%%%%%%%%%%%%
\example\label{_console_fino_ex_}
Let $E_1, E_2$ be elliptic curves.
Consider a Kodaira surface $M$, which is a non-K\"ahler, locally conformally
K\"ahler complex surface, obtained as the total space of a principal holomorphic
$E_1$-bundle  $\pi:\; M \arrow E_2$. Such bundles are classified by the first Chern
class of the fibration $c_1(\pi)$. This class  can be identified with the
$d_2$-differential of the corresponding Leray 
spectral sequence, mapping from $H^1(E_1)$ to $H^2(E_2)$.
By Blanchard's theorem (\cite{_Blanchard_}), $M$ is non-K\"ahler if and 
only if $c_1(\pi)$ is non-zero. The Kodaira surface is an example
of a nilmanifold obtained from the Heisenberg group (\ref{kod_nil}).

Let $\omega_{E_2}$ be a K\"ahler form on $E_2$. 
Then $\pi^*(\omega_{E_2})\in \im d_2$ is exact, giving $\pi^*(\omega_{E_2})= d\xi$.
We chose $\xi$ in such a way that $d(\theta)=0$, where $\theta:=I(\xi)$
(\ref{_Hodge_chasing_Corollary_}). 
Consider a trivial complex line bundle $(L, \nabla_0)$ with connection
defined by $\nabla_0 - \theta$; this line bundle is flat, hence equipped with 
a holomorphic structure operator $\bar\6:=\nabla_0^{0,1}- \theta^{0,1}$
Using the standard Hermitian metric, we express its
Chern connection as $\nabla_0 +\theta^{1,0} - \theta^{0, 1}$,
hence the curvature of this bundle satisfies $\Theta_L= -\1 \pi^*(\omega_{E_2})$.
Rescaling $\theta$ and $\omega_{E_2}$ if necessarily, we can assume
that the cohomology class of $\omega_{E_2}$ is integer.
Then $L$ is a pullback of an ample bundle $B$ on $E_2$. 
Now, $H^0(M,L)= H^0(E_2, B)\neq 0$,
but this space can be interpreted as $H^0(\Lambda^{0,*}(M), \bar\6_\theta)$.
In this example the twisted Dolbeault cohomology of $M$ is  non-zero,
but the corresponding twisted Dolbeault cohomology of the Lie algebra vanishes
(\ref{_twi_Dolbeault_Vanish_}).
This gives a counterexample to Console-Fino theorem for twisted cohomology.

%%%%%%%%%%%%%%%%%%%%%%%%%%%%%%%%%%%%%%%%%%%%%%%%%%%%%%%%%%%%

\section{Classification of LCK nilmanifolds}

%%%%%%%%%%%%%%%%%%%%%%%%%%%%%%%%%%%%%%%%%%%%%%%%%%%%%%%%%%%%

\subsection{Ugarte's conjecture}

The main question concerning nilmanifolds with LCK structure is to decide if  the following conjecture is true or not.

\hfill

\conjecture {\bf (\cite{uga})}\label{conj_uga} 
Let $M$ be a differentiable nilmanifold admitting an LCK structure
(not necessarily locally $G$-invariant, see \ref{_LCK_nilma_Definition_}).
Then $M$ is conformally biholomorphic to a quotient $S^1\times H_{2n-1}/\Lambda$ 
with the Vaisman structure described above.

%%%%%%%%%%%%%%%%%%%%%%%%%%%%%%%%%%%%%%%%%%%%%%%%%%%%%%%%%%%%%%%%%%%%%%%%
\subsection{Classification theorem for LCK nilmanifolds}
%%%%%%%%%%%%%%%%%%%%%%%%%%%%%%%%%%%%%%%%%%%%%%%%%%%%%%%%%%%%%%%%%%%%%%%%

For LCK nilmanifolds, Ugarte's conjecture was proven
by H. Sawai in the following form:

\hfill

\theorem {\bf (\cite{saw1})}\label{inv_nil_vai} 
Let
$(M,I)$  be a complex nilmanifold, $M=G/\Lambda$. 
If $(M, I)$ admits an LCK structure, then it is 
conformally equivalent to Vaisman. 
Moreover, it is biholomorphic
to a quotient of the product $(H_{2n-1}\times\R, I)$.

\hfill

We give a new proof of this result, different from 
the one in \cite{saw1}.

\hfill

\pstep
%The original proof consists in identifying the Lie algebra
%of the covering group as being isomorphic with $
%\goth{h}_{2n-1}\times\R$, and this is done by identifying
%a basis with the same brackets. 
We  replace the LCK metric by a locally $G$-invariant
LCK metric. To this
purpose, we apply the averaging trick that
Belgun used for
the Inoue surface $S^+$ (\cite[Proof of Theorem 7]{_Belgun_});
later, it was generalized
by Fino and Grantcharov (\cite{_Fino_Gra_}). This approach
does not work, in general, for solvmanifolds, but it works 
well for all nilmanifolds. For LCK structures on nilmanifolds, this construction
is due to L. Ugarte:

\hfill

%%%%%%%%%%%%%%%%%%%%%%%%%%%%%%%%%%%%%%%%%%%%%%%%%%%%%%%%%%%%
\theorem \label{bel_av} 
{\bf (\cite[Proposition 34]{uga})}\\
Let $M=G/\Lambda$ be a complex nilmanifold 
admitting an LCK structure $(\omega,\theta)$. Then it
also admits a structure of LCK nilmanifold. In other
words, there exists a left-invariant LCK structure
on $G$ which induces an LCK structure on $M$.

\hfill

\proof As $G$ admits a co-compact lattice, it is
unimodular, hence it admits a bi-invariant measure. Let
$d\mu$ be a bi-invariant volume element on $M$, and
suppose $\vol(M)=1$. Consider a left-invariant 1-form 
form $\check \theta$ on $\g =\Lie(G)$ such that the corresponding
1-form $\theta_0$ on $G/\Lambda$ is cohomologous to $\theta$.
Such a 1-form exists because 
$H^1(G/\Lambda, \R)= H^1(\Lambda^*(\g^*), d)$
(\cite{nomizu}).
Replacing $\omega$ by a conformally equivalent
LCK form $\omega_0$, we can assume that
$d\omega_0=\omega_0\wedge \theta_0$.
Let $D\subset G$ be the fundamental domain of the action
of $\Lambda$.
Given left-invariant vector fields $X, Y$ on $G$, the 2-form
$$\check \omega(X,Y):=\int_D \omega_0(X,Y)d\mu$$
defines an Hermitian structure on the Lie algebra $\g$ of $G$. Moreover, as 
\begin{multline*}
d(\check \omega)(X, Y, Z)\\=  
-\int_D \omega_0([X,Y],Z)d\mu - \int_D
\omega_0(X,[Y,Z])d\mu + \int_D \omega_0(Y,[X,Z])d\mu \\= - \int_D
d(\omega_0)(X,Y,Z) d\mu  = - \int_D \theta_0\wedge\omega_0
(X,Y,Z) d\mu= -\check \omega\wedge\check\theta(X,Y,Z),
\end{multline*}
it follows that $\check \omega$ is an LCK form on $\g$. \endproof

\hfill

{\bf Proof of \ref{inv_nil_vai}, Step 2:}
Now we prove that any LCK nilmanifold is Vaisman.
Using \ref{_Hodge_chasing_Corollary_}, we obtain that $\omega=d_\theta(\tau)$,
where  $d_\theta(I\tau)=0$.

Applying Dixmier and Alaniya's theorem again,
we obtain that $d_\theta(I\tau)=0$ implies $I\tau=d_\theta(v)$,
where $v\in \Lambda^0(\g^*)$ is a constant. Therefore,
$\omega= d_\theta d^c_\theta(\const)$, which is the
equation for the LCK manifold with potential.\footnote
{Note that the constant $\const$
	should be positive by \ref{_strictly_negative_pots_Theorem_}.}
However, as shown in \ref{pot_gau}, an LCK manifold with
potential and constant $|\theta|$ is Vaisman.

\hfill

{\bf   Step 3:} 
We show that the Lee field
$\theta^\sharp\in \g$ and the anti-Lee field
$I(\theta^\sharp)\in \g$ generate an ideal in $\g$. 
This observation also follows from 
\cite[Theorem 3.12]{_Fino_Gra_Ve_}.

Let $\Sigma=\langle \theta^\sharp,I(\theta^\sharp)\rangle$
be the canonical foliation on $M$,
and ${\goth s}\subset \g$ the corresponding subspace of
the Lie algebra of $G$. Since $\theta^\sharp$ and $I(\theta^\sharp)$
are Killing and have constant length, their trajectories
are geodesics (\ref{_canon_foli_totally_geodesic_Remark_}).
Let $S\subset \Iso(M)$ be the group of isometries of $M$
generated by exponentials of $s\in \goth s$.

By construction, $S$ is a subgroup of $G$ which acts on
$M= G/\Lambda$. Therefore, it is invariant under conjugation with
$\Lambda$, and the Lie algebra $\goth s \subset \g$ of $S$ is invariant under the
adjoint action of $\Lambda$ on $\g$. 
Since $G$ is the Mal\v cev completion of $\Lambda$, its image in $\End(\g)$ 
coincides with the Zariski closure of the image of $G$ under the
adjoint action (see Property 1.5 for the Mal\v cev 
completion functor in \cite{_Grunewald_O'Halloran_}). Therefore, $\goth s$ is $G$-invariant.

\hfill

{\bf   Step 4:} 
Now we prove that for any Vaisman nilmanifold $M=G/\Lambda$,
the group $G$ is the product of $\R$ with the Heisenberg
group. Consider the ideal $\goth s\subset \g$ associated with the
canonical foliation $\Sigma$ (Step 3). The leaf space of $\Sigma$
is K\"ahler because $M= G/\Lambda$ is Vaisman.
However, this leaf space is a nilmanifold
associated with the Lie algebra $\g/\goth s$.
A K\"ahler nilpotent Lie algebra is abelian
(\cite{bens_gor}), hence $\g$ is a central
extension of an abelian algebra. Since
$M$ is Vaisman, it is locally isometric to a product
of $\R$ and a Sasakian manifold, hence
the corresponding Lie algebra is a product
of $\R$ and a 1-dimensional central extension\index{extension!central}
of $\R^{2n-2}$. This finishes the proof.
\endproof

\hfill

\remark The general case of \ref{conj_uga} was proven only for
nilmanifolds of Vaisman type, by
G. Bazzoni, \cite{baz}.

\hfill

{\bf Acknowledgment:} Liviu Ornea thanks IMPA for financial support and excellent research environment during the preparation of this work. Both authors thank Anna Fino for carefully reading a first draft of the paper, Sergiu Moroianu for observing an error in a first proof of \ref{_twi_Dolbeault_Vanish_}, and the referees for their very useful suggestions concerning both the content and the references of the paper.

\hfill

{\small

\noindent {\sc Liviu Ornea\\
University of Bucharest, Faculty of Mathematics, \\14
Academiei str., 70109 Bucharest, Romania}, and:\\
{\sc Institute of Mathematics "Simion Stoilow" of the Romanian
Academy,\\
21, Calea Grivitei Str.
010702-Bucharest, Romania\\
\tt lornea@fmi.unibuc.ro, \ \  liviu.ornea@imar.ro}

\hfill

\noindent {\sc Misha Verbitsky\\
{\sc Instituto Nacional de Matem\'atica Pura e
              Aplicada (IMPA) \\ Estrada Dona Castorina, 110\\
Jardim Bot\^anico, CEP 22460-320\\
Rio de Janeiro, RJ - Brasil }\\
also:\\
Laboratory of Algebraic Geometry, \\
Faculty of Mathematics, National Research University 
Higher School of Economics, Department of Mathematics, 9 Usacheva str., Moscow, Russia.
}\\
\tt verbit@verbit.ru, verbit@impa.br }

\end{document}